\documentclass[11pt,thmsa]{article}
\usepackage{amsmath}
\usepackage{amssymb}
\usepackage{graphicx}
\usepackage{epstopdf}

\input{tcilatex}
\begin{document}

\author{F. Carbonell\thanks{%
Biospective Inc., Montreal, Canada}, Y. Iturria-Medina\thanks{%
Montreal Neurological Institute, Canada}, J.C. Jimenez\thanks{%
Instituto de Cibernetica, Matematica y Fisica, La Habana, Cuba}}
\title{Multiple shooting-Local Linearization method for the identification
of dynamical systems}
\date{2015}
\maketitle

\begin{abstract}
The combination of the multiple shooting strategy with the generalized
Gauss-Newton algorithm turns out in a recognized method for estimating
parameters in ordinary differential equations (ODEs) from noisy discrete
observations. A key issue for an efficient implementation of this method is
the accurate integration of the ODE and the evaluation of the derivatives
involved in the optimization algorithm. In this paper, we study the
feasibility of the Local Linearization (LL) approach for the simultaneous
numerical integration of the ODE and the evaluation of such derivatives.
This integration approach results in a stable method for the accurate
approximation of the derivatives with no more computational cost than the
that involved in the integration of the ODE. The numerical simulations show
that the proposed Multiple Shooting-Local Linearization method recovers the
true parameters value under different scenarios of noisy data.
\end{abstract}


\section{Introduction}

Ordinary differential equations (ODEs) are extensively used for modeling the
temporal evolution of complex dynamical systems in dissimilar fields such as
physics, economy, ecology, biology, chemistry and social sciences \cite%
{Chicone2006}. Typically, these ODEs contain parameters that are associated
to phenomenological factors that control the basic variables interplay of
the models. However, the values of such parameters are usually unknown and
must be determined in such a way that the models reproduce the observed
experimental data at best. Despite a time series analysis of observed
experimental data can determine useful quantities that characterize the
system dynamics (e.g., Lyapunov exponents, attractor dimension), identifying
the system structure and estimating the corresponding parameters would be a
matter of greater practical value. Thus, an accurate estimation of the non
observed states and models's parameters is not only critical to reproduce
and describe a given dynamic behavior but also to understand the underlying
causes of the analyzed processes. This is of particular importance for ODEs
describing chaotic dynamics, where the trajectories of interest are very
sensitive to small perturbations of the parameters and initial values (\cite%
{Baake1992}, \cite{Kallrath1993}, \cite{Voss2004}, \cite{Abarbanel2009}). In
this circumstance, a major challenge is to find a proper numerical
integrator able to preserve the stability of the solutions in situations of
parameter-dependent instabilities in such a way that allows an accurate
estimation of these parameters from noisy chaotic observations.

Several strategies have been proposed for dealing with the parameter
estimation problem in ODEs given a set of noisy observations. Among them,
the so-called Initial Value approach is perhaps the most known. In this
approach, the estimated parameters are those that minimize the least square
errors resulting from fitting the numerical solution of the corresponding
Initial Value problem to the given observation data. However, as it has been
pointed out in \cite{Bock1981}, \cite{Bock1983}, \cite{Cao2011}, the
estimators resulting from this approach are very sensitive to the initial
guess of the parameters and usually turn out only local optimum solutions. A
class of estimation methods that overcome this drawback was originally
introduced in \cite{Bock1981} and it is currently known as the Boundary
Value approach (see, e.g., \cite{Bock1983}, \cite{Kallrath1993}, \cite%
{Leineweber2003a}, \cite{Peifer2007}). This approach has two distinctive
components: 1) the introduction of several multiple shooting nodes for
solving the ODE as multiple Initial Value Problems (IVPs) in smaller
subintervals, and 2) the solution of a constrained least squares problem in
an augmented set of parameters. The main advantage of this multiple shooting
strategy is that the whole observed data can be easily used to bring
information about the true solution of the ODE \cite{Bock1983}. Thus, the
solution of the multiple IVP remains close to the true solution since the
initial iteration of the optimization algorithm. In this way, \ the
influence of the poor initial parameter estimates is considerably reduced.
Besides, the splitting of the integration interval into multiple
subintervals limits the error propagation and allows parameter estimation
even for chaotic systems (\cite{Baake1992}, \cite{Kallrath1993}). Despite
the introduction of additional variables seems to yield a more complicated
estimation procedure, it is actually increasing computational efficiency and
numerical stability of the estimation method \cite{Bock1983}, \cite%
{Kallrath1993}. A third estimation strategy, called nonparametric, employs
nonparametric functions to represent the unknown solutions of the ODEs
(see,e.g., \cite{Varah1982}, \cite{Ramsay2007}, \cite{Brunel2008}, \cite%
{Cao2011}, \cite{Wu2012}). Typically, this class of estimators require two
levels of optimization. The lower level approximates non parametric
functions to the ODE trajectories conditional on the ODE parameters, while
the upper optimization level does the estimation of the parameters of
interest. Clearly, as compared to the previous two approaches, this
procedure increases the computational burden of the parameters estimation
process.

As remarked in \cite{Cao2011}, a common difficulty of all these estimation
strategies is the numerical computation of the derivatives of the
trajectories with respect to the parameters of the ODE. With this respect,
three main approximations have been commonly employed. The simplest one,
finite differences, also called external differentiation \cite{Bock1981},
\cite{Bock1983} is not usually recommended due to the high computational
cost required for achieving numerically stable derivatives (see further
discussion in \cite{Peifer2007}). The second one, called internal
differentiation, consists on differentiating the numerical integrator
corresponding to the original differential equation \cite{Bock1981}, \cite%
{Bock1983}, \cite{HairerNorsettWanner200803}. In general, internal
differentiation is a mechanism less computationally intensive than the
external differentiation but, it might introduce also high computational
cost in the case of implicit integrators or integrators defined trough some
numerical derivatives. The third approach (\cite{Bock1981}, \cite{Bock1983},
\cite{Peifer2007}) consists on approximating the variational equations that
describe the temporal evolution of the required derivatives, which must be
integrated simultaneously to original equation. As in the second kind of
approximation, this can be also computationally intensive for certain types
of numerical schemes.

In this paper, we study the feasibility of the Local Linearization (LL)
technique (see, e.g., \cite{jimenez02AMC}, \cite{jimenez05AMC}) for the
simultaneous numerical integration of the IVPs and the evaluation of the
numerical derivatives that appear in the multiple shooting method. In
previous works \cite{Pedroso2003}, \cite{Donnet2007}, \cite{Ginart2010} this
LL technique has been successfully applied for the parameter estimation of
ODEs in the context of the Initial-Value approach. This has been possible
thanks to the convenient trade-off between the numerical accuracy, stability
and computational cost of the LL integrators and their capability of
preserve a number of dynamical behaviors of the ODEs, which became relevant
for the parameter estimation. In addition to this and following the ideas
used in \cite{carbonell02AMC} for the computation of the Lyapunov Exponents,
the LL technique can be used for the numerical integration of the
variational equations associated to the derivative with respect to the
parameters and initial conditions with no more computational cost than the
that involved in the integration of the ODE. Therefore, the application of
the LL technique for identification of ODEs in the framework of Boundary
Value approach is also attractive.

The paper is organized as follows. In Section \ref{MS_Section}, the
essentials on the Multiple Shooting strategy and the generalized
Gauss-Newton algorithm are presented. Section \ref{MS-LL_Section} is focused
in the link of the LL technique to the multiple shooting method. The
resulting algorithm for the parameter estimation is summarized in this
section as well. The performance of the Multiple Shooting-Local
Linearization method is presented in Section \ref{Numerical_Section}
throughout three numerical examples. Finally, some discussion and
conclusions are presented in the last two sections.

\section{Multiple Shooting Method\label{MS_Section}}

Let us consider the $d$-dimensional ODE%
\begin{equation}
\overset{.}{\mathbf{x}}=\mathbf{f}(t,\mathbf{x,p});\ t\in \lbrack t_{0},T]
\label{EqOr}
\end{equation}%
depending on a $p$-dimensional vector $\mathbf{p}$ of parameters, where $%
\mathbf{f}:\mathbb{R}\times \mathbb{R}^{d}\times \mathbb{R}^{p}\rightarrow
\mathbb{R}^{d}$ is a smooth function.

A typical estimation problem for ODEs consist of finding optimal values for
the parameters $\mathbf{p}$ based on the observation of some values of the
state variable $\mathbf{x}$ contaminated with noise (i.e., data points).
That is, suppose that a number of $N$ observed data points $\mathbf{z}_{i}$
related to the state variables $\mathbf{x}$ and parameters $\mathbf{p}$ via
the observation equation
\begin{equation}
\mathbf{z}_{i}=\mathbf{g}(t_{i}^{\ast },\mathbf{x}(t_{i}^{\ast }),\mathbf{p}%
)+\epsilon _{i},  \label{ObsEq}
\end{equation}%
are given at the time instants $t_{i}^{\ast }\in \lbrack t_{0},T]$, $%
i=1,...,N$, where $\mathbf{g}:\mathbb{R}\times \mathbb{R}^{d}\times \mathbb{R%
}^{p}\rightarrow \mathbb{R}^{v}$ is a smooth function, and $\epsilon _{i}$
denotes the measurement errors. If the measurement errors are assumed
independent, Gaussian distributed with zero mean and known variance $\sigma
^{2}$, then the minimization of the weighted least-squares objective function%
\begin{equation*}
J(\mathbf{p})=\sum\limits_{i=1}^{N}\sum\limits_{j=1}^{v}\sigma ^{-2}(%
\mathbf{z}_{i}^{j}-\mathbf{g}^{j}(t_{i}^{\ast },\mathbf{x}(t_{i}^{\ast },%
\mathbf{p}),\mathbf{p}))^{2}
\end{equation*}%
with respect to $\mathbf{p}$ yields a maximum likelihood estimator for the
parameters of the ODE (\ref{EqOr}).

\subsection{Nonlinear optimization problem}

Formally, the least squares problem described so far is a unconstrained
optimization problem of the type%
\begin{equation*}
\min_{\mathbf{p}}\{\left\Vert \mathbf{F}_{1}(\mathbf{p})\right\Vert
_{2}^{2}\},
\end{equation*}%
where $\mathbf{F}_{1}(\mathbf{p})=vec(\mathbf{M}(\mathbf{p}))$ is a $Nv$%
-dimensional vector, $\mathbf{M}(\mathbf{p})$ is a $N\times v$ matrix with
entries $\mathbf{M}^{ji}(\mathbf{p})$ $=$ $\sigma ^{-1}(\mathbf{z}_{i}^{j}-%
\mathbf{g}^{j}(t_{i}^{\ast },\mathbf{x}(t_{i}^{\ast }),\mathbf{p}))$ for all
$i=1,...,N$ and $\ j=1,...,v$, and $vec(.)$ denotes the vectorization
operator.

However, in many applications, certain initial/boundary problems as those
that appear in control engineering problems (see, e.g., \cite{Hartl1995})
additional requirements for the solutions and parameters must be satisfied.
Mathematically, these restrictions are represented by a vector of (component
wise) equality and/or inequality conditions of the form%
\begin{equation*}
\mathbf{R}(t_{1}^{\ast },\mathbf{x}(t_{1}^{\ast }),...,t_{N}^{\ast },\mathbf{%
x(}t_{N}^{\ast }\mathbf{)},\mathbf{p})=0\text{ or }\geq 0.
\end{equation*}%
In this situation, our estimation problem is reformulated as a constrained
optimization problem of the form%
\begin{equation}
\min_{\mathbf{p}}\{\left\Vert \mathbf{F}_{1}(\mathbf{p})\right\Vert _{2}^{2}%
\text{ }|\text{ }\mathbf{R}_{2}(\mathbf{p})=0;\text{ }\mathbf{R}_{3}(\mathbf{%
p})\geq 0\},  \label{OptimizationProblem}
\end{equation}%
for certain functions $\mathbf{R}_{2}$ and $\mathbf{R}_{3}$.

The multiple shooting approach for solving the optimization problem (\ref%
{OptimizationProblem}) consists on the introduction of $m+1$ grid points $%
t_{0}=\tau _{0}<...<\tau _{m}=T$ on the interval $[t_{0},T]$ and new
parameters $\mathbf{s}_{k}=\mathbf{x}(\tau _{k})$, $k=0,...,m$ such that the
solution of the original equation (\ref{EqOr}) can be approximated by the
solution of a set of independent initial value problems%
\begin{eqnarray}
\overset{.}{\mathbf{x}} &=&\mathbf{f}(t,\mathbf{x,p});\ t\in \lbrack \tau
_{k},\tau _{k+1}]  \label{IVP} \\
\mathbf{x}(\tau _{k}) &=&\mathbf{s}_{k}\text{.}  \notag
\end{eqnarray}%
which, in principle, generate a discontinuous trajectory $\{\mathbf{x}%
(t;\tau _{k},\mathbf{s}_{k},\mathbf{p}),$ $t\in \lbrack \tau _{k},\tau
_{k+1}),$ $k=0,...,m-1\}$. These introduced shooting values $\mathbf{s}_{k}$
act as new parameters for the associated optimization problem (\ref%
{OptimizationProblem}) that should be solved for the augmented parameters $%
\mathbf{q=(p,s}_{0},...,\mathbf{s}_{m}\mathbf{)}$. Thus, the optimization
problem (\ref{OptimizationProblem}) is rewritten as%
\begin{equation}
\min_{\mathbf{q}}\{\left\Vert \mathbf{F}_{1}(\mathbf{q})\right\Vert _{2}^{2}%
\text{ }|\text{ }\mathbf{F}_{2}(\mathbf{q})=0;\text{ }\mathbf{R}_{3}(\mathbf{%
q})\geq 0\},  \label{ExtOptimizationProblem}
\end{equation}%
where the vector-valued function $\mathbf{F}_{2}$ contains the equality
restrictions $\mathbf{R}_{2}$ and the continuity conditions

\begin{equation}
\mathbf{c}_{k}=\mathbf{x}(\tau _{k+1};\tau _{k},\mathbf{s}_{k},\mathbf{p})-%
\mathbf{s}_{k+1}=\mathbf{0,}\text{ }k=0,...,m-1.  \label{ContinuityCond}
\end{equation}%
Notice that, the purpose of the imposed continuity conditions (\ref%
{ContinuityCond}) is to guarantee the continuity of the final approximated
solution of the original equation (\ref{EqOr}) rather than updating the
shooting values $\mathbf{s}_{k}$ from interval to interval in (\ref{IVP}).
In fact, the initial value problems (\ref{IVP}) can be independently solved
following a proper parallel running implementation.

\subsection{Linearized optimization problem}

Clearly, (\ref{ExtOptimizationProblem}) represents a very large constrained
non-linear optimization problem that need to be solved via iterative
methods. As originally proposed in \cite{Bock1981}, the damped generalized
Gauss-Newton method is a suitable choice. Thus, starting with initial guess $%
\mathbf{q}^{(0)}=(\mathbf{p}^{(0)},\mathbf{s}_{0}^{(0)},...,\mathbf{s}%
_{m}^{(0)}),$ the Gauss-Newton iteration is given by%
\begin{equation}
\mathbf{q}^{(l+1)}=\mathbf{q}^{(l)}+\alpha _{l}\Delta \mathbf{q}_{l},\text{ }%
l=0,1,...,  \label{GaussNewton}
\end{equation}%
where $0<\alpha _{l}\leq 1$ is a local damping parameter, the increment $%
\Delta \mathbf{q}_{l}$ is the solution of the linearized problem%
\begin{equation*}
\min_{\Delta \mathbf{q}_{l}}\left\{ \left\Vert \mathbf{F}_{1}(\mathbf{q}%
^{(l)})+\frac{\partial \mathbf{F}_{1}}{\partial \mathbf{q}}(\mathbf{q}%
^{(l)})\Delta \mathbf{q}_{l}\right\Vert _{2}^{2}\text{ }|\text{ }\mathbf{F}%
_{2}(\mathbf{q}^{(l)})+\frac{\partial \mathbf{F}_{2}}{\partial \mathbf{q}}(%
\mathbf{q}^{(l)})\Delta \mathbf{q}_{l}=0;\text{ }\mathbf{R}_{3}(\mathbf{q}%
^{(l)})+\frac{\partial \mathbf{R}_{3}}{\partial \mathbf{q}}(\mathbf{q}%
^{(l)})\Delta \mathbf{q}_{l}\geq 0\right\} ,
\end{equation*}%
and the iteration stops when the absolute error condition%
\begin{equation*}
\left\Vert \mathbf{q}^{(l+1)}-\mathbf{q}^{(l)}\right\Vert \leq \varepsilon
\end{equation*}%
holds for certain given tolerance $\varepsilon >0$.

As pointed out in \cite{Bock1981}, it is convenient to choose some of the
observation time points $t_{i}^{\ast }$ as member of the set of multiple
shooting grid points $\tau _{0}<...<\tau _{m}.$ Thus, the choice of the
initial parameters $\mathbf{(s}_{0}^{(0)},...,\mathbf{s}_{m}^{(0)}\mathbf{)}$
can be based on the prior information given by the observation data points,
which is a recognized advantage of the multiple shooting approach. In that
way, despite being discontinuous, the initial trajectory $\{\mathbf{x}%
(t;\tau _{k},\mathbf{s}_{k}^{(0)},\mathbf{p}^{(0)}),$ $t\in \lbrack \tau
_{k},\tau _{k+1})$, $k=0,...,m-1\}$ can remains relatively close to the
observed data points.

For simplicity in our exposition, from now on we will confine to the
equality constrained case. However, as pointed out in
(\cite{Bock1981}), the following results can be straightforwardly
extended to the inequality
constrained case. Thus, the optimal solution of the linearized problem%
\begin{equation}
\min_{\Delta \mathbf{q}_{l}}\left\{ \left\Vert \mathbf{F}_{1}(\mathbf{q}%
^{(l)})+\frac{\partial \mathbf{F}_{1}}{\partial \mathbf{q}}(\mathbf{q}%
^{(l)})\Delta \mathbf{q}_{l}\right\Vert _{2}^{2}\text{ }|\text{ }\mathbf{F}%
_{2}(\mathbf{q}^{(l)})+\frac{\partial \mathbf{F}_{2}}{\partial \mathbf{q}}(%
\mathbf{q}^{(l)})\Delta \mathbf{q}_{l}=0\right\}
\label{LinearOptimizationProblem}
\end{equation}%
is given by%
\begin{equation}
\Delta \mathbf{q}_{l}=-\left( \frac{\partial \mathbf{F}}{\partial \mathbf{q}}%
(\mathbf{q}^{(l)})\right) ^{+}\mathbf{F}(\mathbf{q}^{(l)}),
\label{Increment_Direction}
\end{equation}%
where $\mathbf{F}=\left(
\begin{array}{c}
\mathbf{F}_{1} \\
\mathbf{F}_{2}%
\end{array}%
\right) $,
\begin{equation}
\frac{\partial \mathbf{F}}{\partial \mathbf{q}}=\left(
\begin{array}{ccccc}
\frac{\partial \mathbf{F}_{1}}{\partial \mathbf{s}_{0}} & \frac{\partial
\mathbf{F}_{1}}{\partial \mathbf{s}_{1}} & \cdots & \frac{\partial \mathbf{F}%
_{1}}{\partial \mathbf{s}_{m}} & \frac{\partial \mathbf{F}_{1}}{\partial
\mathbf{p}} \\
\frac{\partial \mathbf{R}_{2}}{\partial \mathbf{s}_{0}} & \frac{\partial
\mathbf{R}_{2}}{\partial \mathbf{s}_{1}} & \cdots & \frac{\partial \mathbf{R}%
_{2}}{\partial \mathbf{s}_{m}} & \frac{\partial \mathbf{R}_{2}}{\partial
\mathbf{p}} \\
\frac{\partial \mathbf{c}_{0}}{\partial \mathbf{s}_{0}} & -\mathbf{I}_{d} &
\cdots & \mathbf{0} & \frac{\partial \mathbf{c}_{0}}{\partial \mathbf{p}} \\
\vdots & \ddots & \ddots & \vdots & \vdots \\
\mathbf{0} & \cdots & \frac{\partial \mathbf{c}_{m-1}}{\partial \mathbf{s}%
_{m-1}} & -\mathbf{I}_{d} & \frac{\partial \mathbf{c}_{m-1}}{\partial
\mathbf{p}}%
\end{array}%
\right)  \label{BigJacobian}
\end{equation}%
and
\begin{equation}
\left( \frac{\partial \mathbf{F}}{\partial \mathbf{q}}\right) ^{+}=\left(
\begin{array}{cc}
\mathbf{I} & \text{ \ }\mathbf{0}%
\end{array}%
\right) \left(
\begin{array}{cc}
\left( \frac{\partial \mathbf{F}_{1}}{\partial \mathbf{q}}\right) ^{\prime }%
\frac{\partial \mathbf{F}_{1}}{\partial \mathbf{q}} & \text{ \ }\left( \frac{%
\partial \mathbf{F}_{2}}{\partial \mathbf{q}}\right) ^{\prime } \\
\frac{\partial \mathbf{F}_{2}}{\partial \mathbf{q}} & \mathbf{0}%
\end{array}%
\right) ^{-1}\left(
\begin{array}{cc}
\left( \frac{\partial \mathbf{F}_{1}}{\partial \mathbf{q}}\right) ^{\prime }
& \text{ \ }\mathbf{0} \\
\mathbf{0} & \text{ \ }\mathbf{I}%
\end{array}%
\right)  \label{GeneralizedInv}
\end{equation}%
denotes a generalized inverse of the Jacobian $\frac{\partial \mathbf{F}}{%
\partial \mathbf{q}}$ (i.e. $\left( \frac{\partial \mathbf{F}}{\partial
\mathbf{q}}\right) ^{+}\frac{\partial \mathbf{F}}{\partial \mathbf{q}}\left(
\frac{\partial \mathbf{F}}{\partial \mathbf{q}}\right) ^{+}=\left( \frac{%
\partial \mathbf{F}}{\partial \mathbf{q}}\right) ^{+}$).

\subsection{ Equivalent condensed problem}

A major challenge in the computation of the optimal solution of the
linearized problem (\ref{LinearOptimizationProblem}) is the algebraic
manipulation of the Jacobian matrix (\ref{BigJacobian}) and, in turn, the
computation of the generalized inverse (\ref{GeneralizedInv}). Notice that
the Jacobian (\ref{BigJacobian}) is a high dimensional matrix of dimension
at least $Nv+d(m-1)$, which makes the direct evaluation of the formula (\ref%
{GeneralizedInv}) computationally unfeasible for a large number of either
observed data points or multiple shooting nodes. However, the sparse
structure in the bottom side in the Jacobian (\ref{BigJacobian}) allows a
convenient recursive elimination of the variables $\Delta \mathbf{s}%
_{m},...,\Delta \mathbf{s}_{1}$. Following \cite{Bock1981}, a backward
recursion can be implemented as%
\begin{eqnarray}
\mathbf{U}_{1}^{(m)} &:&=\mathbf{F}_{1},\text{ }\mathbf{P}_{1}^{(m)}:=\frac{%
\partial \mathbf{F}_{1}}{\partial \mathbf{p}},\text{ }\mathbf{S}_{1}^{(m)}:=%
\frac{\partial \mathbf{F}_{1}}{\partial \mathbf{s}_{m}}
\label{Backward_Iter} \\
\mathbf{U}_{2}^{(m)} &:&=\mathbf{R}_{2},\text{ }\mathbf{P}_{2}^{(m)}:=\frac{%
\partial \mathbf{R}_{2}}{\partial \mathbf{p}},\text{ }\mathbf{S}_{2}^{(m)}:=%
\frac{\partial \mathbf{R}_{2}}{\partial \mathbf{s}_{m}}  \notag \\
\text{For }i &=&m,m-1,...,1\text{:}  \notag \\
\mathbf{U}_{1}^{(i-1)} &:&=\mathbf{U}_{1}^{(i)}+\mathbf{S}_{1}^{(i)}\mathbf{c%
}_{i-1},\text{ }\mathbf{P}_{1}^{(i-1)}:=\mathbf{P}_{1}^{(i)}+\mathbf{S}%
_{1}^{(i)}\left( \frac{\partial \mathbf{c}_{i-1}}{\partial \mathbf{p}}%
\right) ,\text{ }\mathbf{S}_{1}^{(i-1)}:=\frac{\partial \mathbf{F}_{1}}{%
\partial \mathbf{s}_{i-1}}+\mathbf{S}_{1}^{(i)}\left( \frac{\partial \mathbf{%
c}_{i-1}}{\partial \mathbf{s}_{i-1}}\right)  \notag \\
\mathbf{U}_{2}^{(i-1)} &:&=\mathbf{U}_{2}^{(i)}+\mathbf{S}_{2}^{(i)}\mathbf{c%
}_{i-1},\text{ }\mathbf{P}_{2}^{(i-1)}:=\mathbf{P}_{2}^{(i)}+\mathbf{S}%
_{2}^{(i)}\left( \frac{\partial \mathbf{c}_{i-1}}{\partial \mathbf{p}}%
\right) ,\text{ }\mathbf{S}_{2}^{(i-1)}:=\frac{\partial \mathbf{R}_{2}}{%
\partial \mathbf{s}_{i-1}}+\mathbf{S}_{2}^{(i)}\left( \frac{\partial \mathbf{%
c}_{i-1}}{\partial \mathbf{s}_{i-1}}\right) ,  \notag
\end{eqnarray}%
which transforms the problem (\ref{LinearOptimizationProblem}) into the
equivalent condensed problem%
\begin{equation}
\min_{\Delta \mathbf{s}_{0},\Delta \mathbf{p}}\left\{ \left\Vert \mathbf{U}%
_{1}^{(0)}+\mathbf{S}_{1}^{(0)}\Delta \mathbf{s}_{0}+\mathbf{P}%
_{1}^{(0)}\Delta \mathbf{p}\right\Vert _{2}^{2}\text{ }|\text{ }\mathbf{U}%
_{2}^{(0)}+\mathbf{S}_{2}^{(0)}\Delta \mathbf{s}_{0}+\mathbf{P}%
_{2}^{(0)}\Delta \mathbf{p}=0\right\}  \label{CondensedProblem}
\end{equation}%
in the variables $\Delta \mathbf{s}_{0}$ and $\Delta \mathbf{p}$. Notice
that, as compared to (\ref{LinearOptimizationProblem}), the condensed
problem (\ref{CondensedProblem}) is of lower dimension due to the
elimination of the variables $\Delta \mathbf{s}_{m},...,\Delta \mathbf{s}%
_{1} $.

Depending on the nature of the original optimization problem (\ref%
{OptimizationProblem}), the solution of the condensed problem can be
simplified in several ways. The simplest situation is the one where no
equality constrains are required. In this case, the solution of the
condensed problem can be found by solving the system of normal equations $(%
\mathbf{X}^{\prime }\mathbf{X)\beta }=\mathbf{X}^{\prime }\mathbf{y}$ with $%
\mathbf{X}=\left(
\begin{array}{cc}
\mathbf{S}_{1}^{(0)} & \mathbf{P}_{1}^{(0)}%
\end{array}%
\right) ,$ $\mathbf{y}=-\mathbf{U}_{1}^{(0)}$ and $\mathbf{\beta =}\left(
\begin{array}{c}
\Delta \mathbf{s}_{0} \\
\Delta \mathbf{p}%
\end{array}%
\right) $. Another simple situation is where there are no equality
constrains other than an initial condition $\mathbf{x}(t_{0})=\mathbf{x}_{0}$
for the equation (\ref{EqOr}). In this case, it can be seen that $\mathbf{s}%
_{0}=\mathbf{x}_{0}$, $\Delta \mathbf{s}_{0}=\mathbf{0}$ and the condensed
problem is solved with $\mathbf{X}=\mathbf{P}_{1}^{(0)}$ and $\mathbf{\beta =%
}\Delta \mathbf{p}$. For a more general case of equality constrains, the
condensed problem can be solved by using algorithms specifically designed
for linear least squares problems with linear constrains (see \cite%
{Stoer1971}, \cite{Hanson1982}, \cite{Hanson1986} for instance). Once the
condensed problem has been solved for $\Delta \mathbf{s}_{0}$ and $\Delta
\mathbf{p}$, the remaining variables $\Delta \mathbf{s}_{m},...,\Delta
\mathbf{s}_{1}$ can be obtained by the forward recursion%
\begin{equation}
\Delta \mathbf{s}_{i+1}=\left( \frac{\partial \mathbf{c}_{i}}{\partial
\mathbf{s}_{i}}\right) \Delta \mathbf{s}_{i}+\left( \frac{\partial \mathbf{c}%
_{i}}{\partial \mathbf{p}}\right) \Delta \mathbf{p}+\mathbf{c}%
_{i},i=0,...,m-1.  \label{Forward_Iter}
\end{equation}

\subsection{Damping parameter estimation}

It is well-known that the Gauss-Newton iteration (\ref{GaussNewton}) with $%
\alpha _{l}\equiv 1$ guarantees local convergence to a solution $\mathbf{q}%
^{\ast }$ of the problem (\ref{ExtOptimizationProblem}). However, in
practical applications, it is not possible to choose initial parameters
guess $\mathbf{q}^{(0)}$ for guaranteeing iteration convergence to the
optimal global solution. Thus, in order to extend the global convergence
domain, the damping parameter $0<\alpha _{l}\leq 1$ should be chosen to
unsure the decreasing of an appropriate level function $L(\mathbf{q})$
(i.e., $L(\mathbf{q}^{(l+1)})<L(\mathbf{q}^{(l)})$). As pointed out in \cite%
{Bock1981}, this monotonicity test is only feasible when the increment $%
\Delta \mathbf{q}_{l}$ is a descent direction of the level function at $%
\mathbf{q}^{(l)}$. An appropriate choice of $L$ is then given by the locally
defined natural level functions (\cite{Deuflhard1974}, \cite{Bock1981})%
\begin{equation*}
L_{l}(\mathbf{q})=\frac{1}{2}\left\Vert \left( \frac{\partial \mathbf{F}}{%
\partial \mathbf{q}}(\mathbf{q}^{(l)})\right) ^{+}\mathbf{F}(\mathbf{q}%
)\right\Vert _{2}^{2},
\end{equation*}%
for which $\frac{\partial L_{l}}{\partial \mathbf{q}}(\mathbf{q}%
^{(l)})=-\Delta \mathbf{q}_{l}$ (i.e., $\Delta \mathbf{q}_{l}$ is the
steepest descent direction of $L_{l}$ at $\mathbf{q}^{(l)}$).

The damping parameter $\alpha _{l}$ is then determined by%
\begin{equation}
\min_{\alpha _{l}}\{\mathbf{L}_{l}(\mathbf{q}^{(l)}+\alpha _{l}\Delta
\mathbf{q}_{l})\},  \label{LineSearch}
\end{equation}%
which can be solved by any line search algorithm (\cite{Al-Baali1986}, \cite%
{More1994}, \cite{Zhang2004}). Notice that a line search algorithm is also
an iterative procedure that would require additional evaluation of the
function $\mathbf{F}$ at some points $\mathbf{q}^{(l)}+\alpha
_{l}^{(u)}\Delta \mathbf{q}_{l}$, $u=0,1,...$. Correspondingly, an extra
computationally burden appears during the numerical evaluation of the terms
\begin{equation*}
\mathbf{r}_{l}^{u}=\left( \frac{\partial \mathbf{F}}{\partial \mathbf{q}}(%
\mathbf{q}^{(l)})\right) ^{+}\mathbf{F}(\mathbf{q}^{(l)}+\alpha
_{l}^{(u)}\Delta \mathbf{q}_{l}).
\end{equation*}%
Indeed, by using similar arguments to the ones employed for deriving (\ref%
{CondensedProblem}), we can easily observe that the term $-\mathbf{r}%
_{l}^{u} $ is the optimal solution of the linear squares problem%
\begin{equation*}
\min_{\mathbf{r}_{l}^{u}}\left\{ \left\Vert \mathbf{F}_{1}(\mathbf{q}%
^{(l)}+\alpha _{l}^{(u)}\Delta \mathbf{q}_{l})+\frac{\partial \mathbf{F}_{1}%
}{\partial \mathbf{q}}(\mathbf{q}^{(l)})\mathbf{r}_{l}^{u}\right\Vert
_{2}^{2}\text{ }|\text{ }\mathbf{F}_{2}(\mathbf{q}^{(l)}+\alpha
_{l}^{(u)}\Delta \mathbf{q}_{l})+\frac{\partial \mathbf{F}_{2}}{\partial
\mathbf{q}}(\mathbf{q}^{(l)})\mathbf{r}_{l}^{u}=0\right\} ,
\end{equation*}%
which can be solved by reducing it to a corresponding condensed problem.

In order to avoid the intensive evaluations required in full line search
algorithms, we have employed a modified line search method (\cite{Bock1981},
\cite{Peifer2007}) that naturally adapts to the geometry of the problem.
Specifically, the modified line search algorithm consists on finding an
upper bound for the natural level function $L_{l}(\mathbf{q})$, evaluated at
$\mathbf{q=q}^{(l)}+\alpha _{l}\Delta \mathbf{q}_{l}$, which is given by
(see details in \cite{Bock1981} and \cite{Peifer2007})%
\begin{equation*}
L_{l}(\mathbf{q}^{(l)}+\alpha _{l}\Delta \mathbf{q}_{l})\leq \left( 1-\alpha
_{l}+\alpha _{l}^{2}w(\mathbf{q}^{(l)},\alpha _{l})\right) ^{2}L_{l}(\mathbf{%
q}^{(l)}),
\end{equation*}%
where $w(\mathbf{q},\alpha )$ is a function that characterizes the
nonlinearity of the optimization problem (\ref{LineSearch}). The importance
of $w(\mathbf{q},\alpha )$ is given by the fact that (see proof in \cite%
{Peifer2007}), for an arbitrarily chosen $\eta \in (0,2]$, any $\alpha
_{l}\in (0,\alpha ^{\ast }]$ satisfies the required descending property
\begin{equation}
L_{l}(\mathbf{q}^{(l)}+\alpha _{l}\Delta \mathbf{q}_{l})\leq L_{l}(\mathbf{q}%
^{(l)}),  \label{Descending_Prop}
\end{equation}%
where $\alpha ^{\ast }$ is given by%
\begin{equation}
\alpha ^{\ast }=\min \left( 1,\frac{\eta }{w(\mathbf{q}^{(l)},\alpha ^{\ast
})\left\Vert \Delta \mathbf{q}_{l}\right\Vert }\right) .  \label{damp_est}
\end{equation}%
Since $w$ is unknown a priori, an estimator is given by%
\begin{equation}
w(\mathbf{q}^{(l)},\alpha _{l})=2\frac{\left\Vert \left( \frac{\partial
\mathbf{F}}{\partial \mathbf{q}}(\mathbf{q}^{(l)})\right) ^{+}\mathbf{F}(%
\mathbf{q}^{(l)}+\alpha _{l}\Delta \mathbf{q}_{l})-(1-\alpha _{l})\Delta
\mathbf{q}_{l}\right\Vert }{\left\Vert \alpha _{l}\Delta \mathbf{q}%
_{l}\right\Vert ^{2}}.  \label{w_est}
\end{equation}%
Then, a predictor-corrector procedure can be constructed from the two
previous expressions. That is, starting with an estimate $w(\mathbf{q}%
^{(l-1)},\alpha _{l-1})$ from the previous Gauss-Newton iteration $l-1$, the
initial guess $\alpha _{l}^{(0)}$ is determined according to%
\begin{equation*}
\alpha _{l}^{(0)}=\min \left( 1,\frac{\eta }{w(\mathbf{q}^{(l-1)},\alpha
_{l-1})\left\Vert \Delta \mathbf{q}_{l}\right\Vert }\right) .
\end{equation*}%
If the descending property (\ref{Descending_Prop}) holds, then we should
take $\alpha _{l}=\alpha _{l}^{(0)}$ as the optimal damping parameter.
Otherwise, $w$ has to be re-estimated from (\ref{w_est}) with $\alpha
_{l}=\alpha _{l}^{(0)}$ and the process has to be repeated until the
descending property (\ref{Descending_Prop}) be satisfied (see a detailed
algorithm implementation in \cite{Peifer2007}).

\section{Multiple Shooting - Local Linearization method\label{MS-LL_Section}}

Since analytical solutions $\mathbf{x}$ of the ODE (\ref{EqOr}) are
generally unknown, the objective function $J(\mathbf{p})$ is typically
approximated by%
\begin{equation}
\widetilde{J}(\mathbf{p})=\sum\limits_{i=1}^{N}\sum\limits_{j=1}^{v}\sigma
^{-2}(\mathbf{z}_{i}^{j}-\mathbf{g}^{j}(t_{i}^{\ast },\widetilde{\mathbf{x}}%
(t_{i}^{\ast },\mathbf{p}),\mathbf{p}))^{2},  \label{Approx. Obj. Function}
\end{equation}%
where $\widetilde{\mathbf{x}}(t_{i}^{\ast },\mathbf{p})$ denotes a numerical
approximation to $\mathbf{x}(t_{i}^{\ast })$. Therefore, numerical
approximations to functions $\mathbf{F}_{1}$ and $\mathbf{F}_{2}$ as well as
theirs derivatives are needed for evaluating the iteration (\ref{GaussNewton}%
). In this section it is shown how the initial value problems (\ref{IVP}) as
well as the corresponding variational equations respecting to the initial
value and the parameters are numerically approximated by the so-called Local
Linearization approach.

\subsection{Local Linearization integrators\label{Section LL integrators}}

In addition to the IVP (\ref{IVP}), let us consider the associated
variational problems corresponding to the initial value $\mathbf{s}_{k}$%
\begin{eqnarray}
\overset{.}{\mathbf{X}^{s_{k}}} &=&\frac{\partial \mathbf{f}}{\partial
\mathbf{x}}(t,\mathbf{x,p})\mathbf{X}^{s_{k}},\ t\in \lbrack \tau _{k},\tau
_{k+1}]  \label{Varx} \\
\mathbf{X}^{s_{k}}(\tau _{k}) &=&\mathbf{I}_{d}\text{,}  \notag
\end{eqnarray}%
for all $k=0,...,m-1$, where $\mathbf{X}^{s_{k}}=\frac{\partial \mathbf{x}}{%
\partial \mathbf{s}_{k}}$. Here, by definition, $\mathbf{X}^{s_{k}}\equiv
\mathbf{0}_{d}$ for $t\notin \lbrack \tau _{k},\tau _{k+1}]$. Consider also
the variational problem corresponding to the parameters $\mathbf{p}$%
\begin{eqnarray}
\overset{.}{\mathbf{X}^{p}} &=&\frac{\partial \mathbf{f}}{\partial \mathbf{x}%
}(t,\mathbf{x,p})\mathbf{X}^{p}+\frac{\partial \mathbf{f}}{\partial \mathbf{p%
}}(t,\mathbf{x,p}),\ t\in \lbrack \tau _{k},\tau _{k+1}]  \label{Varp} \\
\mathbf{X}^{p}(\tau _{i}) &=&\mathbf{0}_{d\times p}\text{,}  \notag
\end{eqnarray}%
where $\mathbf{X}^{p}=$ $\frac{\partial \mathbf{x}}{\partial \mathbf{p}}$.

Denote by $\Upsilon ^{k}(h)=\{\tau _{k}\leq t_{n}^{k}\leq \tau
_{k+1}:n=0,1,...,N_{k}\}$ a time discretization of the subinterval
\thinspace $\lbrack \tau _{k},\tau _{k+1}]$ with \thinspace $t_{0}^{k}=\tau
_{k}$, $t_{N_{k}}^{k}=\tau _{k+1}$, $h_{n}^{k}=t_{n+1}^{k}-t_{n}^{k}\leq h$
for $h>0$, and satisfying $t_{i}^{\ast }\in \Upsilon ^{k}(h)$ for those
observation time points $t_{i}^{\ast }$ such that $\tau _{k}\leq t_{i}^{\ast
}\leq \tau _{k+1}$. Since the observation time points $t_{i}^{\ast }$, $%
i=1,...,N$ have a fix location over the interval $\left[ t_{0},T\right] $,
any time discretization $\Upsilon ^{k}(h)$ containing more than 2
observation time points does not likely have equally spaced time points $%
t_{n}^{k}$ over the interval $[\tau _{k},\tau _{k+1}]$. Thus, a numerical
integration with a fix step size $h$ is, usually, unfeasible. Instead, an
adaptive step size strategy is in order. For the remaining of our
exposition, it is assumed that the time discretization $\Upsilon ^{k}(h)$
have been constructed under the adaptive step size strategy proposed in \cite%
{Sotolongo2014} for the LL integrators with relative and absolute tolerances
$RelTol$ and $AbsTol$. A slight modification to this adaptive strategy for
including the fix observation time points $t_{i}^{\ast }\in \Upsilon
^{k}(h)\,$\ has been implemented here.

The Local Linear approximation $\mathbf{y}$ to the solution $\mathbf{x}$ of (%
\ref{IVP}) is obtained from the local (piece-wise) linearization of the
function $\mathbf{f}$ respecting to $\mathbf{x}$ and $t$, and the exact
computation of the resulting linear IVP%
\begin{eqnarray}
\overset{.}{\mathbf{y}} &=&\mathbf{f}(t_{n},\mathbf{y}_{t_{n}^{k}},\mathbf{p}%
)+\frac{\partial \mathbf{f}}{\partial \mathbf{x}}(t_{n}^{k},\mathbf{y}%
_{t_{n}^{k}},\mathbf{p})(\mathbf{y}-\mathbf{y}_{t_{n}^{k}})\text{ }\mathbf{+}%
\frac{\partial \mathbf{f}}{\partial t}(t_{n}^{k},\mathbf{y}_{t_{n}^{k}},%
\mathbf{p})(t-t_{n}^{k}),\text{ }t\in \left[ t_{n}^{k},t_{n+1}^{k}\right]
\label{LL_ODE} \\
\mathbf{y}(t_{n}^{k}) &=&\mathbf{y}_{t_{n}^{k}},  \notag
\end{eqnarray}%
with $\mathbf{y}(t_{0}^{k})=\mathbf{y}_{t_{0}^{k}}=\mathbf{s}_{k}$ for all $%
n=0,...,N_{k}$ (see, e.g., \cite{jimenez02AMC},\cite{jimenez05AMC}).

By following the same ideas used in \cite{carbonell02AMC} for computing the
Lyapunov Exponents, the derivatives $\mathbf{X}^{s_{k}}$ and $\mathbf{X}^{p}$
can be approximated by the solution of the variational equations%
\begin{eqnarray}
\overset{.}{\mathbf{Y}^{s_{k}}} &=&\frac{\partial \mathbf{f}}{\partial
\mathbf{x}}(t_{n}^{k},\mathbf{y}_{t_{n}^{k}},\mathbf{p})\mathbf{Y}^{s_{k}},%
\text{ }t\in \left[ t_{n}^{k},t_{n+1}^{k}\right]  \label{LL_Varx} \\
\mathbf{Y}^{s_{k}}(t_{n}^{k}) &=&\mathbf{Y}_{t_{n}^{k}}^{s_{k}}  \notag
\end{eqnarray}%
and%
\begin{eqnarray}
\overset{.}{\mathbf{Y}^{p}} &=&\frac{\partial \mathbf{f}}{\partial \mathbf{x}%
}(t_{n}^{k},\mathbf{y}_{t_{n}^{k}},\mathbf{p})\mathbf{Y}^{p}+\frac{\partial
\mathbf{f}}{\partial \mathbf{p}}(t_{n}^{k},\mathbf{y}_{t_{n}^{k}},\mathbf{p}%
),\text{ }t\in \left[ t_{n}^{k},t_{n+1}^{k}\right]  \label{LL_Varp} \\
\mathbf{Y}^{p}(t_{n}^{k}) &=&\mathbf{Y}_{t_{n}^{k}}^{p},  \notag
\end{eqnarray}%
respectively, with $\mathbf{Y}^{s_{k}}(t_{0}^{k})=\mathbf{Y}%
_{t_{0}^{k}}^{s_{k}}=\mathbf{I}_{d}$ and $\mathbf{Y}^{p}(t_{0}^{k})=\mathbf{Y%
}_{t_{0}^{k}}^{p}=\mathbf{0}_{d\times p}$. Notice that, by construction, $%
\mathbf{Y}_{t_{n}^{r}}^{s_{k}}\equiv \mathbf{0}_{d}$ for $r\neq k,$ $%
n=0,1,...,N_{r}$.

The solutions $\mathbf{y,}$ $\mathbf{Y}^{s_{k}}$ and
$\mathbf{Y}^{p}$ of the equations (\ref{LL_ODE}), (\ref{LL_Varx})
and (\ref{LL_Varp}) can be straightforwardly derived by using their
integral representations obtained in \cite{jimenez02AMC} and
\cite{carbonell02AMC} combined with the formulas for
computing integrals of exponential matrices proposed in \cite{Carbonell2008b}%
. That is,
\begin{equation}
\mathbf{y}_{t_{n+1}^{k}}=\mathbf{y}_{t_{n}^{k}}+\mathbf{E}_{14}(\mathbf{y}%
_{t_{n}^{k}}),\text{ }n=0,...,N_{k}-1  \label{Ll_Scheme_EqOr}
\end{equation}%
\begin{equation}
\mathbf{Y}_{t_{n+1}^{k}}^{s_{k}}=\mathbf{E}_{11}(\mathbf{y}_{t_{n}^{k}})%
\mathbf{Y}_{t_{n}^{k}}^{s_{k}},\text{ }n=0,...,N_{k}-1
\label{Ll_Scheme_EqVars}
\end{equation}%
and%
\begin{equation}
\mathbf{Y}_{t_{n+1}^{k}}^{p}=\mathbf{E}_{11}(\mathbf{y}_{t_{n}^{k}})\mathbf{Y%
}_{t_{n}^{k}}^{p}+\mathbf{E}_{12}(\mathbf{y}_{t_{n}^{k}}),\text{ }%
n=0,...,N_{k}-1  \label{Ll_Scheme_EqVarp}
\end{equation}%
where the vectors $\mathbf{E}_{14}(\mathbf{y}_{t_{n}^{k}})$, $\mathbf{E}%
_{12}(\mathbf{y}_{t_{n}^{k}})$ and the matrix $\mathbf{E}_{11}(\mathbf{y}%
_{t_{n}^{k}})$ are specific block components of the exponential matrix%
\begin{equation*}
\exp (h_{n}^{k}\mathbf{C})=\left[
\begin{array}{cccc}
\mathbf{E}_{11}(\mathbf{y}_{t_{n}^{k}}) & \mathbf{E}_{12}(\mathbf{y}%
_{t_{n}^{k}}) & \mathbf{E}_{13}(\mathbf{y}_{t_{n}^{k}}) & \mathbf{E}_{14}(%
\mathbf{y}_{t_{n}^{k}}) \\
- & - & - & - \\
- & - & - & - \\
- & - & - & -%
\end{array}%
\right]
\end{equation*}%
with $\mathbf{C\in }\mathbb{R}^{(d+p+2)\times (d+p+2)}$ defined as%
\begin{equation*}
\mathbf{C}=\left[
\begin{array}{cccc}
\frac{\partial \mathbf{f}}{\partial \mathbf{x}}(t_{n}^{k},\mathbf{y}%
_{t_{n}^{k}},\mathbf{p}) & \frac{\partial \mathbf{f}}{\partial \mathbf{p}}%
(t_{n}^{k},\mathbf{y}_{t_{n}^{k}},\mathbf{p}) & \frac{\partial \mathbf{f}}{%
\partial t}(t_{n}^{k},\mathbf{y}_{t_{n}^{k}},\mathbf{p}) & \mathbf{f}%
(t_{n}^{k},\mathbf{y}_{t_{n}^{k}},\mathbf{p}) \\
0 & 0 & 0 & 0 \\
0 & 0 & 0 & 1 \\
0 & 0 & 0 & 0%
\end{array}%
\right] .
\end{equation*}%
It is worth noticing here that the numerical implementation of LL schemes (%
\ref{Ll_Scheme_EqOr}), (\ref{Ll_Scheme_EqVars}), (\ref{Ll_Scheme_EqVarp})
reduce to the use of a convenient algorithm for computing matrix
exponentials, e.g., those based on rational Pad\'{e} approximations \cite%
{GolubLoan199610}, the Schur decomposition \cite{GolubLoan199610} or Krylov
subspace methods \cite{sidje1998esp}. The selection of one of them will
mainly depend on the size and structure of the matrix $\mathbf{C}$. For
instance, for many low dimensional system of equations one could use the
algorithm developed in \cite{vanloan1978cii}, which takes advantage of the
special structure of the matrix $\mathbf{C}$. Whereas, for large systems of
equations, the Krylov subspace methods are strongly recommended.

Notice also that the equations (\ref{LL_ODE}), (\ref{LL_Varx}) and (\ref%
{LL_Varp}) are not the result of applying the standard local linearization
technique simultaneously to the set of equations (\ref{IVP}), (\ref{Varx})
and (\ref{Varp}). Instead, an appropriate local linearization approach has
been chosen in order to decouple the system of equations (\ref{IVP}), (\ref%
{Varx}) and (\ref{Varp}). Indeed, (\ref{LL_ODE}) is the local linear
approximation to equation (\ref{IVP}) but equations (\ref{LL_Varx}) and (\ref%
{LL_Varp}) are suitable linear equations with locally constant coefficients.
Nevertheless, it has been proved in \cite{carbonell02AMC} that%
\begin{equation*}
\underset{t\in \lbrack \tau _{k},\tau _{k+1}]}{\sup }\left\Vert \mathbf{X}%
^{s}(t)-\mathbf{Y}^{s}(t)\right\Vert \leq C_{s}^{k}h,
\end{equation*}%
where the constant $C_{s}^{k}$ does not depend on $h$. Correspondingly,
following similar arguments to the ones employed in Theorem 4 of \cite%
{carbonell02AMC}, it can be also proved that%
\begin{equation*}
\underset{t\in \lbrack \tau _{k},\tau _{k+1}]}{\sup }\left\Vert \mathbf{X}%
^{p}(t)-\mathbf{Y}^{p}(t)\right\Vert \leq C_{p}^{k}h,
\end{equation*}%
for certain constant $C_{p}^{k}$. Thus, despite
\begin{equation*}
\underset{t\in \lbrack \tau _{k},\tau _{k+1}]}{\sup }\left\Vert \mathbf{x}%
(t)-\mathbf{y}_{t}\right\Vert \leq C^{k}h^{2},
\end{equation*}%
for certain constant $C^{k}$ (see proof in \cite{jimenez02AMC}), the system
of equations (\ref{LL_ODE})-(\ref{LL_Varp}) has global order of convergence
equal to 1. In other words, the numerical derivatives $\mathbf{X}^{s}$ and $%
\mathbf{Y}^{s}$ can be approximated with global order of convergence 1 and
no extra computationally cost but the one involved in the implementation of
the local linearization schemes. Remarkably, it has been also avoided the
manipulation of second order derivatives like to ones that would certainly
appear with the employ of internal differentiation in the equation (\ref%
{LL_ODE}). Additionally, under request, the Lyapunov exponents of the ODEs
might be straightforwardly approximated from the solution $\mathbf{Y}^{s}$
by following the algorithm developed in \cite{carbonell02AMC}.

\subsection{Parameters estimation algorithm}

The Multiple Shooting-Local Linearization algorithm for estimating the
unknown parameters $\mathbf{p}$ of the model (\ref{EqOr})-(\ref{ObsEq})
proceeds by inserting the LL approximations of the previous subsection into
the minimization objective function (\ref{Approx. Obj. Function}), namely,%
\begin{equation*}
\widetilde{J}(\mathbf{p})=\sum\limits_{i=1}^{N}\sum\limits_{j=1}^{v}\sigma
^{-2}(\mathbf{z}_{i}^{j}-\mathbf{g}^{j}(t_{i}^{\ast },\mathbf{y}%
_{t_{i}^{\ast }},\mathbf{p}))^{2},
\end{equation*}%
where $\mathbf{y}_{t_{i}^{\ast }}$ denotes the LL\ approximation to $\mathbf{%
x}(t_{i}^{\ast })$, $i=1,...,N$. Correspondingly, the continuity constrains $%
\mathbf{c}_{k}$ and additional equality constrains take the form $\mathbf{c}%
_{k}=\mathbf{y}_{t_{N_{k}}^{k}}-\mathbf{s}_{k+1}$, $k=0,...,m-1$, and $%
\mathbf{R}_{2}=\mathbf{R}_{2}(t_{1}^{\ast },\mathbf{y}_{t_{1}^{\ast
}},...,t_{N}^{\ast },\mathbf{y}_{t_{N}^{\ast }},\mathbf{p})$, respectively.
Analogously, the functions $\mathbf{F}_{1}(\mathbf{p})$, $\mathbf{F}_{2}(%
\mathbf{p})$ and $\mathbf{R}_{3}$ of the Section \ref{MS_Section} must be
redefined in terms of the approximations $\mathbf{y,}$ $\mathbf{Y}^{s_{k}}$
and $\mathbf{Y}^{p}$ to $\mathbf{x,}$ $\mathbf{X}^{s_{k}}$ and $\mathbf{X}%
^{p}$. \ Indeed, from now on, $\mathbf{F}_{1}(\mathbf{p})=vec(\widetilde{%
\mathbf{M}}(\mathbf{p}))$ with $\ \widetilde{\mathbf{M}}^{ji}(\mathbf{p})$ $%
=\sigma ^{-1}(\mathbf{z}_{i}^{j}-\mathbf{g}^{j}(t_{i}^{\ast },\mathbf{y}%
_{t_{i}^{\ast }},\mathbf{p}))$,
\begin{eqnarray}
\frac{\partial \mathbf{F}_{1}}{\partial \mathbf{s}_{k}} &=&[\frac{\partial
\mathbf{g}}{\partial \mathbf{y}}(t_{1}^{\ast },\mathbf{y}_{t_{1}^{\ast }},%
\mathbf{p})\mathbf{Y}_{t_{1}^{\ast }}^{s_{k}};\frac{\partial \mathbf{g}}{%
\partial \mathbf{y}}(t_{2}^{\ast },\mathbf{y}_{t_{2}^{\ast }},\mathbf{p})%
\mathbf{Y}_{t_{2}^{\ast }}^{s_{k}};...;\frac{\partial \mathbf{g}}{\partial
\mathbf{y}}(t_{N}^{\ast },\mathbf{y}_{t_{N}^{\ast }},\mathbf{p})\mathbf{Y}%
_{t_{N}^{\ast }}^{s_{k}}],  \label{Evaluations} \\
\frac{\partial \mathbf{F}_{1}}{\partial \mathbf{p}} &=&[\frac{\partial
\mathbf{g}}{\partial \mathbf{y}}(t_{1}^{\ast },\mathbf{y}_{t_{1}^{\ast }},%
\mathbf{p})\mathbf{Y}_{t_{1}^{\ast }}^{p}+\frac{\partial \mathbf{g}}{%
\partial \mathbf{p}}(t_{1}^{\ast },\mathbf{y}_{t_{1}^{\ast }},\mathbf{p});%
\frac{\partial \mathbf{g}}{\partial \mathbf{y}}(t_{2}^{\ast },\mathbf{y}%
_{t_{2}^{\ast }},\mathbf{p})\mathbf{Y}_{t_{2}^{\ast }}^{p}+\frac{\partial
\mathbf{g}}{\partial \mathbf{p}}(t_{2}^{\ast },\mathbf{y}_{t_{2}^{\ast }},%
\mathbf{p});...;  \notag \\
&&\frac{\partial \mathbf{g}}{\partial \mathbf{y}}(t_{N}^{\ast },\mathbf{y}%
_{t_{N}^{\ast }},\mathbf{p})\mathbf{Y}_{t_{N}^{\ast }}^{p}+\frac{\partial
\mathbf{g}}{\partial \mathbf{p}}(t_{N}^{\ast },\mathbf{y}_{t_{N}^{\ast }},%
\mathbf{p})],  \notag \\
\frac{\partial \mathbf{R}_{2}}{\partial \mathbf{s}_{k}} &=&\sum%
\limits_{i=1}^{N}\frac{\partial \mathbf{R}_{2}}{\partial \mathbf{y}_{i}}%
(t_{1}^{\ast },\mathbf{y}_{t_{1}^{\ast }},...,t_{N}^{\ast },\mathbf{y}%
_{t_{N}^{\ast }},\mathbf{p})\mathbf{Y}_{t_{i}^{\ast }}^{s_{k}}  \notag \\
\frac{\partial \mathbf{R}_{2}}{\partial \mathbf{p}} &=&\sum\limits_{i=1}^{N}%
\frac{\partial \mathbf{R}_{2}}{\partial \mathbf{y}_{i}}(t_{1}^{\ast },%
\mathbf{y}_{t_{1}^{\ast }},...,t_{N}^{\ast },\mathbf{y}_{t_{N}^{\ast }},%
\mathbf{p})\mathbf{Y}_{t_{i}^{\ast }}^{p}+\frac{\partial \mathbf{R}_{2}}{%
\partial \mathbf{p}}(t_{1}^{\ast },\mathbf{y}_{t_{1}^{\ast
}},...,t_{N}^{\ast },\mathbf{y}_{t_{N}^{\ast }},\mathbf{p})  \notag \\
\frac{\partial \mathbf{c}_{k}}{\partial \mathbf{s}_{k}} &=&\mathbf{Y}%
_{_{t_{N_{k}}^{k}}}^{s_{k}},  \notag \\
\frac{\partial \mathbf{c}_{k}}{\partial \mathbf{p}} &=&\mathbf{Y}%
_{_{t_{N_{k}}^{k}}}^{p},  \notag
\end{eqnarray}%
where $[.;.;...;.]$ denotes the algebraic operation of concatenating
matrices with equal number of columns by their rows. Here, $\mathbf{Y}%
_{t_{i}^{\ast }}^{s_{k}}$ and $\mathbf{Y}_{t_{i}^{\ast }}^{p}$ denote the LL
approximations to $\mathbf{X}^{s_{k}}(t_{i}^{\ast })$ and $\mathbf{X}%
^{p}(t_{i}^{\ast }),$respectively$.$

The parameters estimation algorithm is then summarized in the following
steps:

\begin{enumerate}
\item Setting $l=0$ and initial guess $\mathbf{q}^{(0)}=(\mathbf{p}^{(0)},%
\mathbf{s}_{0}^{(0)},...,\mathbf{s}_{m}^{(0)})$ for the parameters and
shooting nodes,

\item With $\mathbf{p}=\mathbf{p}^{(l)}$ and $\mathbf{s}_{k}=\mathbf{s}%
_{k}^{(l)},$ $k=1,...,m$, compute $\mathbf{y}_{t_{i}^{\ast }},\mathbf{Y}%
_{t_{i}^{\ast }}^{s_{k}}$ and $\mathbf{Y}_{t_{i}^{\ast }}^{p}$ as indicated
in Section \ref{Section LL integrators} for all $i=1,...,N$. Then, evaluate
the expressions (\ref{Evaluations}),

\item Compute the increments $\Delta \mathbf{q}_{l}$ in (\ref%
{Increment_Direction}) by either direct evaluation of the Jacobian (\ref%
{BigJacobian}) and the generalized inverse (\ref{GeneralizedInv}) or
evaluating the backward and forward iterations (\ref{Backward_Iter}) and (%
\ref{Forward_Iter}) in the condensed problem,

\item Compute the damping parameter $\alpha _{l}$ by the modified line
search algorithm according to (\ref{damp_est})-(\ref{w_est}),

\item Iterate the Gauss-Newton algorithm $\mathbf{q}^{(l+1)}=\mathbf{q}%
^{(l)}+\alpha _{l}\Delta \mathbf{q}_{l}$,

\item Set $l=l+1$ and repeat steps (2)-(5) until $\left\Vert \mathbf{q}%
^{(l+1)}-\mathbf{q}^{(l)}\right\Vert \leq \varepsilon $ for a given
tolerance $\varepsilon >0.$
\end{enumerate}

\subsection{Variance estimation}

In practical situations, the variance $\sigma ^{2}$ of the observation
errors in (\ref{ObsEq}) is also an unknown parameter that should be
estimated, namely, by extending the parameter $\mathbf{p}$ with the
inclusion of $\sigma $. However, since only the function $\mathbf{F}_{1}$
does depend on $\sigma $, the inclusion of $\sigma $ in the Gauss-Newton
iteration process would unnecessarily increase the dimension of the problem.
An alternative estimation for $\sigma $ is then computed as%
\begin{equation*}
\mathbf{\sigma }^{(l)}=\sqrt{\frac{\sum\limits_{i=1}^{N}\sum%
\limits_{j=1}^{v}(\mathbf{z}_{i}^{j}-\mathbf{g}^{j}(t_{i}^{\ast },\mathbf{y}%
_{t_{i}^{\ast }},\mathbf{p}^{(l)}))^{2}}{Nv-p}},\text{ }l=0,1,...,
\end{equation*}%
Obviously, in this case, the estimated $\widehat{\mathbf{p}}$ is not longer
maximum likelihood estimator.

\section{Numerical Experiments\label{Numerical_Section}}

In this section, the performance of the Multiple Shooting-Local
Linearization approach is illustrated through three numerical examples. The
first example, extensively studied in \cite{Kallrath1993}\thinspace , is a
4-dimensional chaotic system defined by a vector field that is linear
respecting to the unknown parameters. The second example corresponds to the
well-known FitzHugh-Nagumo system, which is defined nonlinearly respecting
to the parameters of interest. The last example correspond to the Rikitake
system \cite{Rikitake1958}, which is known for generating chaotic
trajectories fro certain parameters combination. For the three examples, the
parameters were estimated with a stopping tolerance of $\varepsilon =$ $%
10^{-4}$ and the shooting points were selected within the set of the
observed time points $t_{i}^{\ast }$, $i=1,...,N$, in an approximately
equispaced manner. For each $t_{i}^{\ast }$, $i=1,...,N$, the LL
approximations $\mathbf{y}_{t_{i}^{\ast }},\mathbf{Y}_{t_{i}^{\ast
}}^{s_{k}} $ and $\mathbf{Y}_{t_{i}^{\ast }}^{p}$ were adaptively computed
with relative and absolute tolerances $RelTol=10^{-3}$ and $AbsTol=10^{-6}$.

\textbf{Example 1.} Consider the Henon-Heiles system described by the
4-dimensional ODE (see details in \cite{Kallrath1993}):%
\begin{eqnarray*}
\overset{.}{x}_{1} &=&\mathbf{x}_{3} \\
\overset{.}{x}_{2} &=&x_{4} \\
\overset{.}{x}_{3} &=&-ax_{1}-2x_{1}x_{2} \\
\overset{.}{x}_{4} &=&-bx_{2}-x_{1}^{2}-cx_{2}^{2},
\end{eqnarray*}%
with parameters $\mathbf{p}=(a,b,c)$. The \ "true" trajectory in the
interval $[0,10]$ is shown in Figure 1 for $\mathbf{p}=(1,1,-1)$ and initial
condition $\mathbf{x}_{0}=(0,0,0.3,-0.4)$. This \ "true" trajectory $\mathbf{%
x}$ was generated by the Local Linearization method with a fixed step size
of $h=2^{-12}$. A realization of $N$ random observations $\mathbf{z}_{i}$, $%
\,$is generated by randomly selecting $N$ points $t_{i}^{\ast }$, $%
i=1,...,N, $ $\,$in the interval $[0,10]$ (with uniform distribution) and
adding a Gaussian noise with zero mean and variance $\sigma ^{2}$ to the
value $\mathbf{x}(t_{i}^{\ast })$. That is,
\begin{equation*}
\mathbf{z}_{i}=\mathbf{x}(t_{i}^{\ast })+\sigma \epsilon _{i},\text{ }%
\epsilon _{i}\thicksim N(0,1)\text{, }i=1,...,N,
\end{equation*}%
with $N(0,1)$ denoting the Gaussian normal distribution. A number of 1000 of
such realizations were generated for different values of $\sigma $ and $N$.
These 1000 realizations were arranged into 20 batches of 50 realizations
each, where each batch corresponds to a fix distribution of the observation
time points $t_{i}^{\ast }$, $i=1,...,N$. The distribution of the
observation time points then varies from batch to batch. The goal was to
estimate the parameters $\mathbf{p}$, $\mathbf{x}_{0}$ and $\sigma $ in each
realization. For each realization, the initial parameters guesses were set
at $\mathbf{p}^{(0)}=$ $(9,1,2)$ and $\sigma ^{(0)}=1$, and $m=50$ shooting
nodes were distributed over the interval $[0,10]$.

It should be noticed that the integration of this chaotic system with
initial condition $\mathbf{x}_{0}=(0,0,0.3,-0.4)$ and parameter $\mathbf{p}%
=(9,1,2)$ leads to numerically unstable solutions (i.e. numerical
explosions) after $t=4.4$ even with a very small fixed step size of $%
h=2^{-12}$. This evidences that the classical initial value approach
estimation is not suitable in this scenario. Instead, more sophisticated
methods like the multiple shooting approach presented here seems to be a
proper choice.

The estimated parameters are reported in Table 1 as the average within the
batch (i.e. average across 100 realization of fixed observation time points
distribution) and then average and standard deviation across the 20 batches.
Notice that such a summary should not be confounded with the so-called a
posteriori analysis (see \cite{Bock1981}) that is usually carried out for
statistical inference of the estimated parameters (e.g. variance-covariance
matrix and confidence interval for the estimated parameters).
\begin{eqnarray*}
&&%
\begin{tabular}{lllll}
\hline
& \multicolumn{2}{l}{$\ \ \ \ \ \ \ \ \ \ \ \ \ \ \ \ \ \ \ \ \ \sigma =0.05$%
} & \multicolumn{2}{l}{$\ \ \ \ \ \ \ \ \ \ \ \ \ \ \ \ \sigma =0.1$} \\
\hline
$N$ & $\ \ \ \ \ \ \ \ \ \ \ 100$ & $\ \ \ \ \ \ \ \ \ \ 200$ & $\ \ \ \ \ \
\ \ \ \ \ \ 100$ & $\ \ \ \ \ \ \ \ \ \ \ \ 200$ \\ \hline
$\widehat{a}$ & $\ \ \ 1.0002\pm 0.0006$ & $\ \ 0.9998\pm 0.0007$ & $\ \ $\ $%
1.0008\pm 0.0021$ & $\ \ 1.0003\pm 0.0014$ \\ \hline
$\widehat{b}$ & $\ \ \ 0.9986\pm 0.0017$ & $\ \ 1.0003\pm 0.0020$ & $\ \ $\ $%
0.9973\pm 0.0045$ & $\ \ 0.9997\pm 0.0026$ \\ \hline
$\widehat{c}$ & $-0.9993\pm 0.0042$ & $-0.9991\pm 0.0025$ & $-1.0018\pm
0.0068$ & $-0.9989\pm 0.0061$ \\ \hline
$\widehat{\mathbf{x}_{0}}$ & $%
\begin{array}{c}
-0.0001\pm 0.0014 \\
-0.0004\pm 0.0011 \\
\text{ \ }0.2995\pm 0.0006 \\
-0.4005\pm 0.0008%
\end{array}%
$ & $%
\begin{array}{c}
\text{ \ }0.0001\pm 0.0009 \\
\text{ }0.0001\pm 0.0010 \\
\text{ \ }0.2997\pm 0.0005 \\
-0.3999\pm 0.0005%
\end{array}%
$ & $%
\begin{array}{c}
-0.0016\pm 0.0031 \\
-0.0003\pm 0.0022 \\
\text{ \ }0.2990\pm 0.0015 \\
-0.4003\pm 0.0018%
\end{array}%
$ & $%
\begin{array}{c}
\text{ \ }0.0006\pm 0.0019 \\
-0.0001\pm 0.0013 \\
\text{ \ }0.2998\pm 0.0008 \\
-0.4003\pm 0.0015%
\end{array}%
$ \\ \hline
$\widehat{\sigma }$ & $\ \ \ 0.0496\pm 0.0001$ & $\ \ \ 0.0501\pm 0.0001$ & $%
\ \ \ 0.0992\pm 0.0001$ & $\ \ \ 0.0997\pm 0.0001$ \\ \hline
$N.Iter.$ & $\ \ \ 6.4320\pm 1.0058$ & $\ \ \ 5.0220\pm 0.4527$ & $\ \ \
7.5070\pm 0.7993$ & $\ \ \ 6.8910\pm 0.3061$ \\ \hline
\end{tabular}
\\
&&\text{Table 1. Estimated parameters and number of required Gauss-Newton
iterations (}N.Iter.\text{) } \\
&&\text{corresponding to the Henon-Heiles system.}
\end{eqnarray*}%
Figure 1 shows the true trajectory with initial condition $\mathbf{x}%
_{0}=(0,0,0.3,-0.4)$ and $N=100$ noisy observations corresponding to one
realization with $\sigma =0.1$. This figure also shows the approximated
discontinuous trajectory after the first iteration as well the estimated
optimal trajectory after $l=6$ iterations of the Gauss-Newton method$.\,$%
This optimal trajectory corresponds to the estimated parameters $\widehat{%
\mathbf{x}}_{0}=(-0.0231,-0.0008,0.3055,-0.3899)$, $\widehat{\mathbf{p}}%
=(1.0314,0.9839,-1.0101)$ \ and $\widehat{\sigma }=0.1029$. Notice
that the first iteration produces a discontinuous trajectory due to
the continuity conditions (\ref{ContinuityCond}) are unable to be
satisfied at this stage of the optimization process. However, after
only four iterations, the estimated parameters and shooting nodes
produce an optimal continuous trajectory that is quite close to the
true trajectory of the problem.

\begin{figure}[h]
\centering
\includegraphics[width=6.5in]{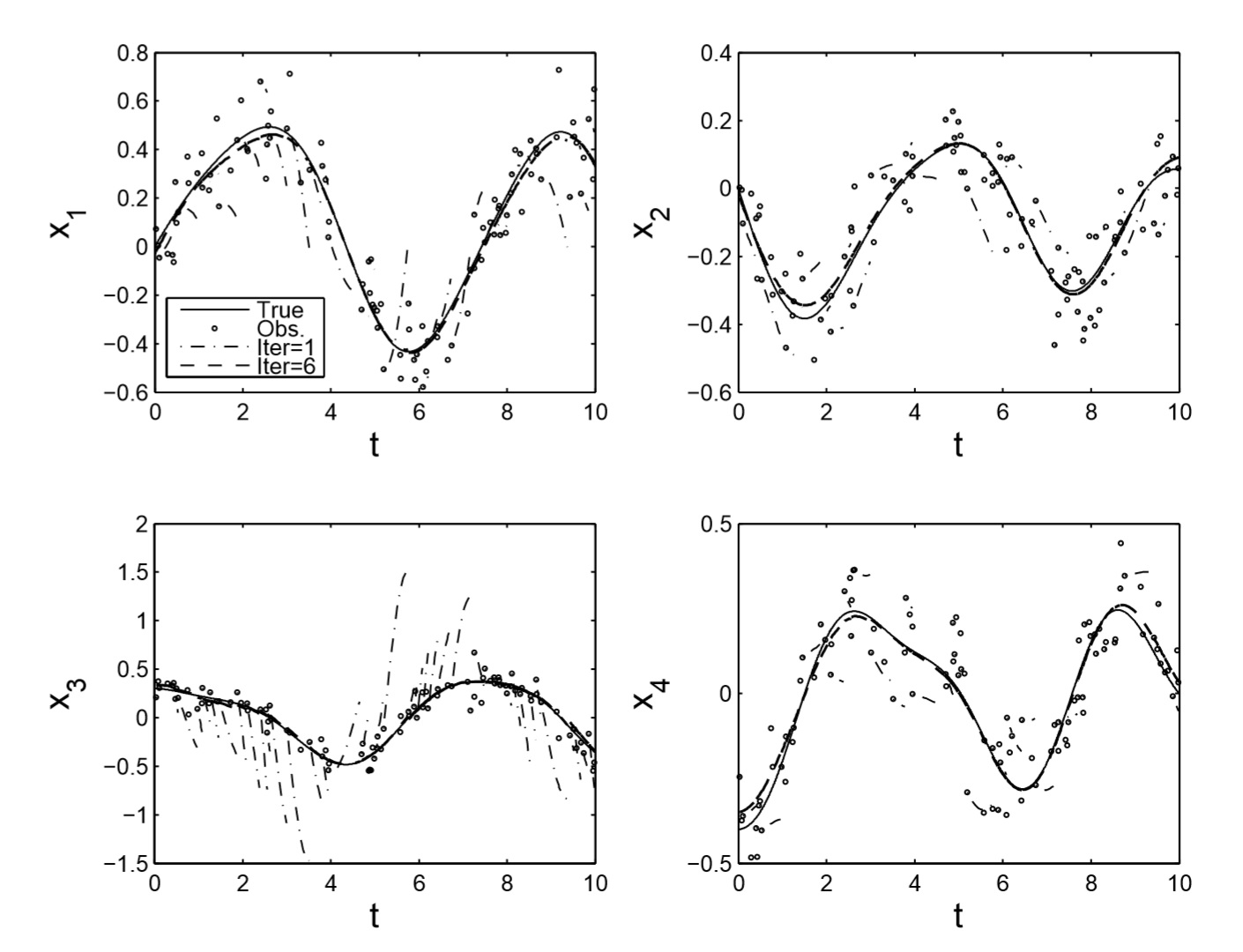}
\caption{Initial and optimal trajectory corresponding to the
Henon-Heiles system with estimated initial condition $
\widehat{\mathbf{x}}_{0}=(-0.0231,-0.0008,0.3055,-0.3899)$ and
parameters $ \widehat{\mathbf{p}}=(1.0314,0.9839,-1.0101)$ \ and
$\widehat{\protect\sigma }=0.1029$.}
\end{figure}

\textbf{Example 2.} Consider the FitzHugh-Nagumo ODE, which is a
simplified version of the well-known Hodgkin--Huxley model for
describing activation
and deactivation dynamics of a spiking neuron:%
\begin{eqnarray*}
\overset{.}{V} &=&c(V-\frac{V^{3}}{3}+R) \\
\overset{.}{R} &=&-\frac{1}{c}(V-a+bR),
\end{eqnarray*}%
where $V$ and $R$ denote the voltage across an axon membrane and the
outwards currents, respectively. Here, $a,b,c$ are parameters to be
estimated from $n=400$ noisy observations of the variable $V$, which were
randomly distributed (with uniform distribution) in the interval $[0,20]$.
Similarly to \cite{Ramsay2007} and \cite{Cao2011}, the true trajectory was
generated with initial values $V(0)=-1$ and $R(0)=1$ and true parameters $\
a=0.2$, $b=0.2$ and $c=3$. The noisy observations were generated by adding a
Gaussian noise with standard deviation $\sigma $ $=0.2$. The initial
parameter guess was set $\mathbf{p}^{(0)}=(a^{(0)},b^{(0)},c^{(0)})=(2,2,5)$
and $\sigma ^{(0)}=1$. A number of $m=50$ shooting nodes were approximately
equispaced over the set of observation time points. Since only the variable $%
V$ is observed in the case and no additional information if available for
the variable $R$ at the shooting points, we set the second component of $%
\mathbf{s}_{k}^{(0)}$ equal to zero for all $k=0,...,m$.

The estimated parameters resulting from 1000 realizations (20 batches of 50
realizations each) were $\widehat{a}=0.2007\pm 0.0023$, $\widehat{b}%
=0.1932\pm 0.0068$, $\widehat{c}=2.9794\pm 0.0113,\widehat{\mathbf{x}}%
_{0}=(-1.0019\pm 0.0123,1.0085\pm 0.0121)$ and $\widehat{\sigma }=0.2016\pm
0.0010$. Figure 2 shows the true, initial and estimated trajectories after $%
l=29$ iterations. Notice that a larger number of iterations were required in
this case probably caused by the very bad (far away from the true
trajectory) initial guess of the second component in the shooting nodes. The
estimated trajectory corresponds to parameters with values $\widehat{a}%
=0.1971$, $\widehat{b}=0.2210$ and $\widehat{c}=2.9716.$

\begin{figure}[h]
\centering
\includegraphics[width=6.5in]{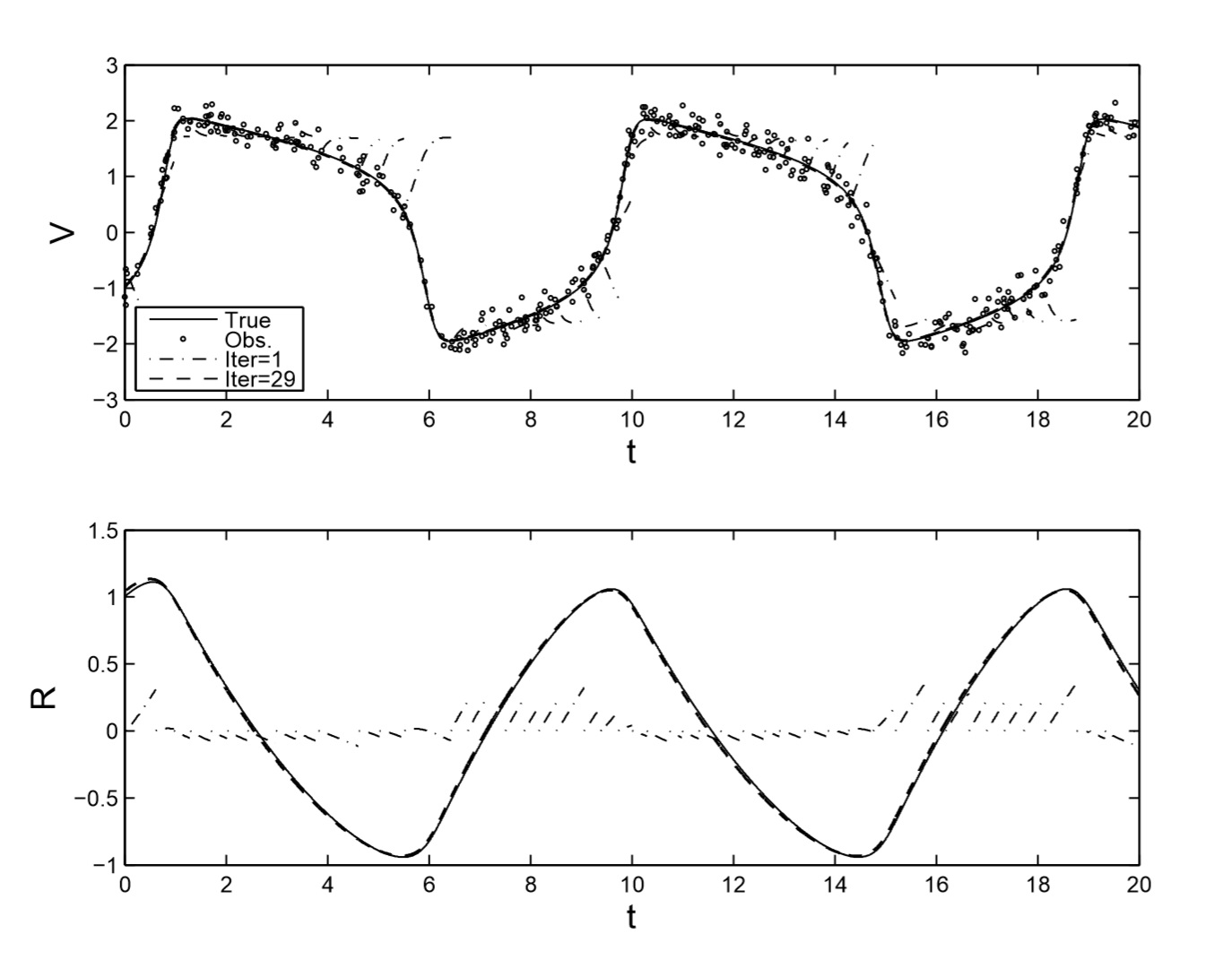}
\caption{Initial and optimal trajectory
corresponding to the FitzHugh-Nagumo system with estimated parameters $%
\widehat{a}=0.1971$, $\widehat{b}=0.2210$ and $\widehat{c}=2.9716.$}
\end{figure}

\textbf{Example 3. }Consider the Rikitake model defined by the ODE%
\begin{eqnarray*}
\overset{.}{x}_{1} &=&-\mu x_{1}+x_{2}x_{3} \\
\overset{.}{x}_{2} &=&-\alpha x_{1}-\mu x_{2}+x_{1}x_{3} \\
\overset{.}{x}_{3} &=&1-x_{1}x_{2},
\end{eqnarray*}%
which was originally introduced by \cite{Rikitake1958} to explain
geomagnetic polarity reversals. The model consists of coupled, self-excited
disc dynamos, where the parameter $\mu >0$ and $\alpha $ represent the
resistive dissipation and the difference of the angular velocities of two
dynamo discs, respectively. Despite the physical meaning of $\mu $ is still
not clear, estimates of geophysically plausible value for $\mu $ vary
between $10^{-3}$ and $10$ \cite{Ito1980}. Most of the studies for
explaining the dynamical behavior of the Rikitake system focus on the
parameter space determined by the pairs \thinspace $(\mu ,K)\,$, where $%
\alpha =\mu (K^{2}-K^{-2})~$(see \cite{Ito1980} for instance). Thus,
combinations of the pairs $(\mu ,K)$ produce different dynamical regimes,
like the chaotic regime determined by $\mu =0.5$ and $\alpha =0.46125$ $%
(K=1.25).$

For this example, the parameters $\mu $ and $\alpha $ are going to be
estimated from $N=200$ noisy observations of the three variable, randomly
distributed (with uniform distribution) in the interval $[0,40]$. A\ "true"
trajectory was simulated with initial value $\mathbf{x}_{0}=(-2,-2,0)$ and
the noisy observations were generated by adding a Gaussian noise with
standard deviation $\sigma $ $=0.1.$ The initial parameter guesses were set
at $\mathbf{p}^{(0)}=(\mu ^{(0)},\alpha ^{(0)})=(5,5)$ and $\sigma ^{(0)}=1$%
. The following table presents the estimated parameters for different
numbers of shooting nodes, including the case $m=0$ corresponding to the
Initial Value approach. The estimated parameters are reported by the average
and standard deviation over 100 different realizations of the observations $%
\mathbf{z}_{i}$ with a fix (random) distribution of the $N$ observation time
points $t_{i}^{\ast }$, $i=1,...,N$. Notice that the average and standard
deviation were calculated only across those realizations where the
estimation algorithm converged after a maximum number of 50 iterations. In
fact, this table also shows the required number of Gauss-Newton iterations ($%
N.Iter.$) that the algorithm needed to converge as well as the percentage of
convergence ($\%Conv.$).
\begin{eqnarray*}
&&%
\begin{tabular}{lcccc}
\hline
$m$ & $60$ & $40$ & $30$ & $20$ \\ \hline
$\widehat{\mu }$ & $\ \ \ 0.5001\pm 0.0005$ & $\ \ 0.5005\pm 0.0137$ & $\ \
0.4998\pm 0.0199$ & $\ \ 0.4965\pm 0.0359$ \\ \hline
$\widehat{\alpha }$ & $\ \ \ 0.4613\pm 0.0010$ & $\ \ 0.4623\pm 0.0194$ & $\
\ 0.4581\pm 0.0248$ & $\ \ 0.4351\pm 0.0719$ \\ \hline
$\widehat{\mathbf{x}_{0}}$ & $%
\begin{array}{c}
-1.9992\pm 0.0269 \\
-1.9993\pm 0.0204 \\
\text{ \ }0.0002\pm 0.0490%
\end{array}%
$ & $%
\begin{array}{c}
-2.0121\pm 0.1066 \\
-2.0551\pm 0.1654 \\
\text{ \ }0.0074\pm 0.0927%
\end{array}%
$ & $%
\begin{array}{c}
-2.0183\pm 0.1583 \\
-2.0297\pm 0.1927 \\
\text{ \ }0.0311\pm 0.1913%
\end{array}%
$ & $%
\begin{array}{c}
-2.1491\pm 0.3754 \\
-2.1228\pm 0.4045 \\
\text{ \ }-0.4141\pm 0.5947%
\end{array}%
$ \\ \hline
$\widehat{\sigma }$ & $\ \ 0.0999\pm 0.0033$ & $\ \ 0.1324\pm 0.1926$ & $\ \
\ 0.1798\pm 0.3397$ & $\ \ \ 0.3764\pm 0.6537$ \\ \hline
$N.Iter.$ & $\ \ \ 15.51\pm 1.13$ & $\ \ \ \ \ 18.73\pm 3.94$ & $\ \ \ \ \
19.07\pm 5.72$ & $\ \ \ \ \ 27.35\pm 7.26$ \\ \hline
$\%Conv.$ & $\ \ \ \ \ \ \ \ \ \ \ \ \ 99$ & $\ \ \ \ \ \ \ \ \ \ \ \ 88$ & $%
\ \ \ \ \ \ \ \ \ \ \ \ 77$ & $\ \ \ \ \ \ \ \ \ \ \ \ 17$ \\ \hline
\end{tabular}
\\
&&\text{Table 2. Estimated parameters, number of required Gauss-Newton
iterations (}N.Iter.\text{), } \\
&&\text{and percentage of convergence (}\%Conv.\text{) corresponding to the
Rikitake system.}
\end{eqnarray*}

Notice that as the number of shooting nodes decreases, the estimated
parameters become less accurate and the number of non convergent
realizations increases. In fact, the simulations corresponding to $m=10$ and
$m=0$ showed no convergent realization at all, which evidences the efficacy
of the multiple shooting method as compared to the Initial Value approach.
Importantly, recall that, due to the equivalent condensed problem,
increasing the number of shooting nodes does not increase the dimensionality
of the optimization problem. Therefore, as a rule of thumb, it is
recommendable to employ the multiple shooting approach with a relatively
large number of shooting nodes, particularly for those system showing
complex, chaotic dynamics.

\section{Discussion}

The methodology presented here can be extended in several ways. As it was
mentioned earlier, only the case of equality constrains for parameters and
state variables has been treated here. For inequality constrains, it is easy
to check that the condensing recursion is exactly the same as for equality
constrains. Therefore, the solution of the condensed problem must be
obtained by more general optimization strategies like active set strategy
(see details in \cite{Bock1981} and \cite{Bock1983}). Additionally, we have
assumed a very simple assumption for the measurements errors that define the
set of observed data points. Namely, uncorrelated and equally distributed
errors have been assumed for the components of the multi-dimensional data.
This scenario can be easily extended to the more general case of correlated
errors by replacing the parameter $\sigma ^{2}$ by a variance-covariance
matrix $\Sigma ~$\ and defining a proper formulation of the function $%
\mathbf{F}_{1}(\mathbf{p})$. Correspondingly, the observed data and
the measurements errors might define more complicated statistical
models like mixed effects models to cover, for instance, the cases
of repeated measures at certain time points and
temporarily-correlated errors.

Finally, the multiple shooting-LL approach can covers a more general class
of models driven by random differential equations (RDE). Essentially, a RDE
is a non autonomous ODE coupled with a stochastic process, which is usually
employed for modelling noisy perturbations of deterministic systems. Thus,
in principle, a RDE can be integrated by applying conventional numerical
methods for ODEs, like the LL integrator presented here \cite{Carbonell2005}%
. In fact, the LL method for RDE has been already successfully applied for
the generation of EEG rhythms by means of realistically coupled neural mass
models \cite{Sotero2007}. A possible extension consists of having more
realistic neural mass models with certain free parameters that could be
estimated from observed EEG data via the multiple shooting approach.

\section{Conclusions}

In this paper we have shown the feasibility of the multiple shooting
approach in combination with local linearization techniques for parameter
estimation in ordinary differential equations. The main advantage of the
proposed approach consists of approximating the numerical derivatives
involved in the multiple shooting scheme by a numerically stable method at
no extra computational burden but the one required for the numerical
integration of the original equations. The performance of the proposed
approach has been evaluated in three different numerical examples. In all
cases, the multiple shooting-local linearization method accurately recovered
the true parameters values.

\bibliographystyle{elsart-num}
\bibliography{acompat,MyPapers,OursPapers,OurUsedBooks,OurUsedPapers}

\end{document}